\documentclass[10pt,twoside]{article}
\usepackage{graphicx}
\usepackage{amsmath}
\usepackage{Latex-document}

\def\volumeno{I}

\title{\bf International Comparisons in \vskip -2mm Mathematics Education:  An Overview}

\author{Gabriele Kaiser\thanks{University of Hamburg, Department of Education,
Von-Melle-Park 8, 20146 Hamburg, Germany. E-mail:
gkaiser@erzwiss.uni-hamburg.de} \quad Frederick K. S.
Leung\thanks{Faculty of Education, University of Hong Kong,
Pokfulam Road, Hong Kong. E-mail: hraslks@hku.hk} \\ Thomas
Romberg\thanks{University of Wisconsin-Madison, Wisconsin Center
for Education Research, 1025 West Johnson Street, Madison,
Wisconsin, USA. E-mail: tromberg@facstaff.wisc.edu} \quad Ivan
Yaschenko\thanks{Moscow Center for Continuous Math Education, B.
Vlas'evskij 11, 121002 Moscow, Russia. E-mail: ivan@mccme.ru}}

\date{\vspace{-8mm}}

\begin{document}

\maketitle

\thispagestyle{first} \setcounter{page}{631}

\begin{abstract}

\vskip 3mm

The paper opens with an overview of the discussion of international comparisons
(including goals) in mathematics education. Afterwards, the two most important
recent international studies, the PISA Study and TIMSS-Repeat, are
described. After a short description of the qualitative-quantitative debate, a
qualitatively oriented small-scale study is described. The paper closes with
reflection on the possibilities and limitations of such studies.

\vskip 4.5mm


\end{abstract}

\vskip 12mm

\section{Goals of comparative studies}

\vskip-5mm \hspace{5mm}

Since the results of the Third International Mathematics and Science Study
(TIMSS) were published in 1996, international comparisons of student performance
in mathematics have gained more and more importance as a consequence of public
and political discussions. The discussions recently have been fueled by the
results published in 2001 of the Programme for International Student Assessment
(PISA). Nevertheless it has to be considered that comparative education has a
long tradition going back to oral reports, as exemplified by Greeks and
Romans. With the beginning of the 19th century, approaches were developed
seeking to identify forces influencing the development of systems of
education. In the 1960s and 1970s, the use of social science methods became
common in order to examine the effect of various factors on educational
development accompanied by a debate on the relative merits of quantitative
versus qualitative studies. We will come back to this discussion later on.

If we look for the goals of comparative education, history shows us that
comparative education serves a variety of goals. It can deepen our understanding
of our own education and society, be of assistance to policymakers and
administrators, and be a valuable component of teacher education
programmes. These contributions can be made through work that is primarily
descriptive as well as through work that seeks to be analytic or explanatory,
through work that is limited to just one or a few nations, and through work that
relies on nonquantitative as well as quantitative data and methods.  Based on
that, Postlethwaite \cite{11} discriminated four major aims of comparative education:
\begin{itemize}
\item
``Identifying what is happening elsewhere that might help improve our own system
of education" (p.xx). Postlethwaite gave several examples, such as the attempt
to identify the principles involved in an innovation like mastery learning
(which has had such success in the Republic of Korea) and grasping the
procedures necessary to implement the mastery principle.

\item ``Describing similarities and differences in educational phenomena between
systems of education and interpreting why these exist" (p.xx). This
comprises the analysis of similarities and differences between systems
of education in goals, in structures, in the scholastic achievement of
age groups and so on, which could reveal important information about the
systems being compared. Studies of these types might describe not only
inputs to and processes within systems but also the philosophy of
systems and outcomes.  The reasons of why certain countries have
particular philosophies and the implications these have in terms of
educational outcomes are areas of both major academic and practical
interest.

\item ``Estimating the relative effects of variables (thought to be
determinants) on outcomes (both within and between systems of education)"
(p.xx). Within education there is a great deal of speculation about what affects
what.  How much evidence, for example, do the people who teach methods at
teacher-training establishments have on the effectiveness of the methods they
promulgate?  What about home versus school effects on outcomes?  These questions
and similar ones are the questions to be dealt with under this perspective.

\item ``Identifying general principles concerning educational effects" (p.xx).
This means that we are aiming at a possible pattern of relationship between
variables within an educational system and an outcome.  In practice, a model
will be postulated whereby certain variables are held constant before we examine
the relationship between other variables and the outcomes.  The resultant
relationship will often be estimated by a regression coefficient. Principles we
detect in an educational system that we analyze that recurs in other systems
might be determined to be a general principle.
\end{itemize}

In mathematics education, there have been a remarkable number of
international comparative studies carried out in the last 30
years.  Robitaille [12, p. 41] believed that the reason for this
might be that---
\begin{quote}
Studies that cross national boundaries provide participating countries with a
broader context within which to examine their own implicit theories, values and
practices.  As well, comparative studies provide an opportunity to examine a
variety of teaching practices, curriculum goals and structures, school
organisational patterns, and other arrangements for education that might not
exist in a single jurisdiction.
\end{quote}
Stigler and Perry \cite[p.~199]{14} emphasized the better understanding of one's own
culture gained through comparative studies:
\begin{quote}
Cross cultural comparison also leads researchers and educators to a more
explicit understanding of their own implicit theories about how children learn
mathematics.  Without comparison, we tend not to question our own traditional
teaching practices, and we may not even be aware of the choices we have made in
constructing the educational process.
\end{quote}

\section{Recent international studies in mathematics}

\vskip-5mm \hspace{5mm}

In the following, we will describe briefly the most important comparative
studies in mathematics education (for details see \cite{4}). Most of the
large-scale studies have been carried out by the International Association for
the Evaluation of Educational Achievement (IEA).

\subsection{From FIMS over SIMS to TIMSS}

\vskip-5mm \hspace{5mm}

The first large-scale international study was the First International
Mathematics Study (FIMS), carried out 1964.  Twelve countries participated in
this study, in which two populations were tested---thirteen-year-olds and
students in the final school year of the secondary school.  In the first
population, the students from Israel, Japan and Belgium received the best
results, and the worst results were achieved by the U.S. students.  In the
second population, a different picture emerged---the youngsters from Israel,
Belgium and England received the best results, and the U.S. students the
worst. Several critics emphasized the important role of the curriculum and
stated that valid comparative results cannot be formulated without considering
curricular aspects.

The second large-scale comparative study in mathematics education was the Second
International Mathematics Study (SIMS), 1980--1982, the results of which were
published at the end of the 1980s and the beginning of the 1990s. Twenty
countries participated in this study, which considered the same age groups as
FIMS and contained a cross-sectional and a longitudinal component.  Considering
the curricular criticisms on FIMS, SIMS discriminated different levels of the
curriculum---the intended curriculum, the implemented curriculum, and the
attained curriculum. In addition, a content by cognitive-behavior grid was
developed, which related the mathematical content with cognitive dimensions such
as computation and comprehension.  On the level of the intended curriculum, the
main results were a significant curricular shift---geometry had lost importance
in contrast to number and geometry. On the level of the implemented curriculum,
the study pointed out the different status of repetition in the different
countries.  On the level of the attained curriculum, the study showed that the
increase in the achievements was remarkable low in many countries.  Gender
differences emerged in many countries, but were not consistent and were smaller
than the differences between the different countries.  SIMS has been criticised
from several perspectives, and even the organisers of SIMS admitted that,
despite the wealth of items, the curricula of many participating countries had
not been covered sufficiently.

The last study in this series is the Third International Mathematics and
Science Study (TIMSS), which was carried out in 1995 in over 40
countries. It examined the achievement of students from three
populations at five grade levels (9-year-olds, 13-year-olds, and
students in the final year of secondary school) in a wide range of
content and performance areas, and it collected contextual information
from students, teachers, and school principals. In considering the
criticisms formulated at SIMS---the unsatisfactory coverage of the
curriculum of the different countries, the focus on quantified outcomes
(the quantified achievement of the students)---TIMSS established several
additional studies:
\begin{itemize}
\item The TIMSS videotape study, which analyzed mathematics lessons in Japan,
Germany, and the United States;

\item The case study project, which collected qualitative information on the
educational systems in Japan, Germany and the United States;

\item The survey of mathematics and science opportunities, a study of
mathematics and science teaching in six countries;

\item The curriculum analysis study, which studied the curricula and textbooks in many countries.
\end{itemize}
With all these additional studies and the high number of countries participating
in the main study, TIMSS remains the largest and most comprehensive study of
educational practice in mathematics and science ever undertaken.

\subsection{TIMSS repeat}

\vskip-5mm \hspace{5mm}

Because of the impact of TIMSS on the international community, it was
decided that a repeat study (TIMSS Repeat or TIMSS-R) be conducted so
that trends in mathematics (and science) achievements could be studied
in an international context.  However, TIMSS was a complicated study
involving testing three populations of students at five grade levels and
in two subject areas.  So it was decided that for TIMSS-R in 1999, only
eighth grade students (i.e., the upper grade of the TIMSS population 2
level (i.e., 13-year-olds)) would be tested.

38 countries participated in TIMSS-R, and of the 38 countries, 26 had
participated in the eighth grade test in TIMSS as well.  So for these 26
countries, comparison between their eighth grade students' results in 1995 and
1999 could be made.  17 of these 26 countries had participated in the fourth
grade test in TIMSS as well, and for these 17 countries, the choice of
replicating the TIMSS study in 1999 means that the students tested in 1999 was
the same cohort of students who took the TIMSS test in 1995.  This thus
constitutes a quasi-longitudinal study and trends in achievement across the
four-year duration can be studied.  And for the 12 countries which did not
participate in TIMSS in 1995, the TIMSS-R results would allow them to compare
their students' achievements with all the TIMSS and TIMSS-R countries.

TIMSS-R, being a repeat study of TIMSS, adopted the TIMSS framework,
which in turn was based on the SIMS framework.  As pointed out above,
SIMS and TIMSS placed special emphasis on the curriculum, and conceived
the curriculum in the three levels of intended, implemented and attained
curricula.  The TIMSS curriculum framework, which was developed through
a lengthy process of negotiation among National Research Coordinators of
all the countries that participated in TIMSS with input from experts in
the field of mathematics education, includes a Content dimension and a
Performance expectations dimension (There is a third dimension known as
Perspectives, but the TIMSS-R international report has not included
results based on analysis of this dimension of the data).  The five
content areas tested included: Fractions and Number Sense (38\% of the
items were devoted to this area); Measurement (15\%); Data
Representation, Analysis and Probability (13\%); Geometry (13\%); and
Algebra (22\%).  The categories under performance expectation were:
Knowing (19\%); Using Routine Procedures (23\%); Using Complex
Procedures (24\%); Investigating and Problem Solving (31\%); and
Communicating and Reasoning (2\%).  The test consisted of items in
multiple-choice and free-response (short answers and extended responses)
formats, and about one quarter of the items and one third of the testing
time were devoted to the free-response items.

Like TIMSS, questionnaires for school principals, teachers and students were
administered to collect data on the variables related to student achievement.
Also, in line with the IEA tradition, rigorous sampling and administration
standards were established and closely monitored by the international study
centre.

Results:

Mathematics Achievement

For the top performing countries in TIMSS-R, the pattern in TIMSS persisted. The
East Asian countries (Singapore, Korea, Chinese Taipei, Hong Kong and Japan)
outperformed their counter-parts in other parts of the world.  Other countries
that achieved well included Belgium (Flemish), Netherlands, Slovak Republic,
Hungary, Canada, Slovenia, Russian Federation, Australia, Finland, Czech
Republic, Malaysia and Bulgaria.  Like TIMSS, one conspicuous finding in TIMSS-R
is the magnitude of the difference in achievements among countries.  The range
of performance among countries is more than three standard deviations, and in
fact, the achievements of all the five top performing East Asian countries were
three standard deviations above that of the lowest scoring country, and those of
the top six countries was more than two standard deviations above those of the
three lowest performing countries.

For most of the countries that participated in both TIMSS and TIMSS-R, there was
not much difference in terms of their relative performance in the two studies.
Only one country had a significant decrease in its performance, and 3 out of the
26 had a significant improvement.

Variables Related to Achievement

System Variables

There was no clear relationship between the wealth of the countries (as measured
by GNP) and their students' achievements.  Although many affluent countries
(Japan, Singapore, Belgium (Flemish), Netherlands, Hong Kong) did very well in
TIMSS-R, some wealthy countries (such as US, Finland, Australia, Italy, and
Israel) did not do as well.  On the other hand, some less affluent countries
(Slovak Republic, Korea and Chinese Taipei) did very well in TIMSS-R.  Nor was
high achievement related to public expenditures on education.  In fact, the
public expenditures on education of the top 9 countries were all less than the
international average percentage of 5.13\%.

For the average over all TIMSS-R countries, students with a higher level of
educational resources at home and in school did better in the TIMSS-R test.
However, we cannot conclude that students from countries with higher levels of
educational resources did better in TIMSS-R.  Some high achieving countries
(such as Singapore and Hong Kong) had relatively few educational resources,
while some countries with high level of resources (such as Israel and US) did
not do very well.  So while there might well be a positive correlation between
educational resources and student achievement within countries, there did not
seem to be a clear relationship between educational resources and achievement
across countries.

Students' Attitudes towards Mathematics

Similar to the finding above, although the TIMSS-R results were consistent with
the findings from the literature that students' positive attitudes towards
mathematics was related with higher achievement within a country, the same
relationship did not hold across countries.  In fact, with the exception of
Singapore, all the top-performing countries had relatively negative attitudes
towards mathematics.

In particular, across countries, a positive self-concept in mathematics did not
seem to be related with higher achievement.  It is noticeable that students from
all the five top-performing East Asian countries had very low self-image of
mathematics.  This suggests that self-image of mathematics or confidence in
doing mathematics may be related to cultural values and is not necessarily
associated with student achievement.

Teacher and Instructional Practices

The same pattern applied to teacher confidence.  Although the TIMSS-R data show
that within a country, teacher confidence was related with student achievement,
teachers from high performing countries did not have particularly high level of
confidence.  In particular, although Japanese students did very well in TIMSS-R,
their teachers had the lowest level of confidence among all the TIMSS-R
countries.

Time devoted to mathematics instruction varied tremendously across countries,
from the lowest of 73 hours and 9\% of the total instructional time, to 222
hours and 17\% of the total instructional time.  But once again, there was no
clear relationship between amount of instructional time and achievement.  In
actual fact, the four countries that devoted most time to mathematics
instruction did badly in TIMSS-R, and the top performing countries were spending
about just half of the time as these countries in instruction.

As far as instructional practices are concerned, teachers from participating
countries reported that the two most predominant activities in their classrooms
were teacher lecture and teacher-guided student practice, which accounted for
nearly half of the class time.  Although solving non-routine problems was
mentioned in the intended curriculum of nearly all countries, teachers in all
TIMSS-R countries, with the exception of Japan, reported that they put
relatively low emphasis on mathematics reasoning and problem solving.

There were other discrepancies between the intended and the implemented
curriculum.  ``For example, curricular goals and aims in 25 countries included
``visualization of three-dimensional shapes" for all or almost all students, but
teachers in only eight countries reported that at least 75 percent of the
students had been taught this topic." (\cite[p.~182]{9})

Lastly, the amount of homework assigned by teachers to students differed
tremendously across countries, but again there was no neat relationship between
the amount of homework and the achievement of students.

Conclusion

As far as the TIMSS-R achievement within countries is concerned, the factors we
discussed above merely confirmed the findings in previous studies, that
achievement is related to educational resources at home and in school, with
students' attitudes towards mathematics, with teachers' confidence etc.  But
what TIMSS-R failed to inform us is how to account for differences in
performance across countries.  As pointed out above, most of the variables that
explained variation of achievement within countries failed to explain the
variation across countries.  This failure perhaps points to the limitation of
questionnaires in getting at the factors for explaining student achievement.
This problem is particularly acute in international studies, where the same term
in a questionnaire may mean very different things in different culture.  There
is also the issue of confoundment between the cultural values and the measuring
instrument.  For example, the finding on negative self-concept in mathematics of
East Asian students above may be due to the stress in the East Asian cultures of
the virtue of humility or modesty.  Children from these countries are taught
from when they are young that one should not be boastful.  This may inhibit
students from rating themselves too highly on the question of whether they think
they do well in mathematics, and so the scores may represent less than what
students are really thinking about themselves.  On the other hand, one's
confidence and self image are something that is reinforced by one's learned
values, and if students are constantly taught to rate themselves low, they may
internalize the idea and may result in really low confidence.

For variables on instructional practices, the TIMSS-R teacher questionnaire
results again did not give us a lot of clues on the instructional practices that
will lead to high achievement.  Probably, it is hard for instructional practices
to be captured by self-reporting questionnaire, and that is why associated
studies such as the TIMSS-R Video Study are so important in this regard.
Video-taping offers a form of cross-cultural documentation which is both true to
the original classroom and amenable to rigorous analysis \cite{15}, and is hence a
much better methodology for studying instructional practices.

What TIMSS-R does tell us is that there exist vast differences in mathematics
achievement across a large number of countries.  Hopefully the realization of
the differences will fuel a search for the factors that contribute to high
achievement rather than a race to top the league table.

\subsection{OECD/PISA}

\vskip-5mm \hspace{5mm}

Developed jointly by the Organisation for Economic Cooperation and Development
(OECD) member countries, the Program for International Student Assessment (PISA)
is designed to monitor, on a regular basis, the literacy of students in reading,
mathematics, and science as they approach the end of secondary school. The first
PISA assessment took place in 2000 with the emphasis on reading, with initial
assessments in both mathematics and science. The next assessment will be in 2003
with the emphasis on mathematics. PISA has been implemented through an
international consortium led by the Australian Council for Educational Research
(ACER) (for details see \cite{10}).

The OECD/PISA Mathematical Literacy Study is concerned with the capacities of
students to analyse, reason, and communicate ideas effectively as they pose,
formulate, solve, and interpret solutions to mathematical problems in a variety
of domains and situations. By focusing on real world problems, the PISA
assessment does not limit itself to the kinds of situations and problems
typically encountered in school classrooms.  In real world settings, few people
are faced with the drill-type of problems that typically appear in school
textbooks and classrooms. Instead, citizens in every country are currently being
bombarded with information on issues such as ``global warming and the greenhouse
effect," ``population growth," ``oil slicks and the seas," ``the disappearing
countryside," and so on, and a relevant question is whether citizens can make
sense of claims and counterclaims on such issues. Of interest for the OECD/PISA
study is whether 15-year-olds (the age when many students are completing their
formal compulsory mathematics learning) can use the mathematics they have been
taught to help make sense of these kinds of issues.

The concept of mathematical literacy, which underlies the OECD/PISA study is defined as---
\begin{quote}
an individual's capacity to identify and understand the role that mathematics
plays in the world, to make well-founded judgements, and to engage in
mathematics, in ways that meet the needs of that individual's life as a
constructive, concerned, and reflective citizen.
\end{quote}

This definition of mathematical literacy is consistent with the broad and
integrative theory about the structure and use of language as reflected in
recent sociocultural literacy studies. The term ``literacy" refers to the human
use of language. In fact, each human language and each human use of language has
both an intricate design tied in complex ways to a variety of functions.  For a
person to be literate in a language implies that the person knows many of the
design resources of the language and is able to use those resources for several
different social functions. Analogously considering mathematics as a language
implies that students not only must learn the design features involved in
mathematical discourse (the terms, facts, signs and symbols, procedures, and
skills in performing certain operations in specific mathematical subdomains and
the structure of those ideas in each subdomain), they also must learn to use
such ideas to solve nonroutine problems in a variety of situations defined in
terms of social functions (making sense of some phenomena). Note that the design
features for mathematics are more than knowing the basic terms, procedures, and
concepts that one is commonly taught in schools. It involves how these features
are structured and used. Unfortunately, one can know a good deal about the
design features of mathematics without knowing either their structure or how to
use those features to solve problems.

PISA assess mathematical literacy in three dimensions:
\begin{enumerate}
\item Content in terms of broad mathematics domains such as chance, change and
growth, space and shape, uncertainty, among others. For PISA 2000, the
mathematics assessment focused on two domains: change and growth, and space and
shape.

\item Three ``competency classes." The OECD/PISA mathematical literacy items
have been developed to assess three classes of student mathematical competency:

Class 1 competencies include those most commonly used on standardised
assessments and classroom tests.  These competencies are knowledge of facts,
common problem representations, recognition of equivalents, recalling of
familiar mathematical objects and properties, performance of routine procedures,
application of standard algorithms and technical skills, manipulation of
expressions containing symbols and formulae in standard form, and computations.

Class 2 competencies include those related to students' planning for problem
solving by drawing connections between the different mathematical content
strands, or from different Big Ideas.  They also include students' abilities to
combine and integrate information in order to tackle and solve ``standard"
problems. Class 2 competencies reflect students' abilities to choose and develop
strategies, to choose mathematical tools, to use multiple methods or steps in
the mathematisation and modelling process. These competencies also reflect
students' abilities to interpret and reflect on the meaning of a solution and
the validity of their work.  Problems that reflect student competencies in Class
2 require students to use appropriate elements from different mathematical
content areas, or from different Big Ideas, in combination with conceptual
thinking and reasoning based on material that does not call for large extensions
of where the student has been before.

Class 3 competencies include those related to students' ability to plan solution
strategies and implement them in problem settings that are more complex and
``original" (or unfamiliar) than those in Class 2. These competencies require
students not only to mathematise more complex problems, but also to develop
original solution models.  Items measuring Class 3 competencies should reflect
students' ability to analyse, to interpret, to reflect on, and to present
mathematical generalisations, arguments and proofs.

\item Situations in which mathematics is used. For PISA, each item was set in
one of five situation-types: personal, educational, occupational, public, and
scientific. The items selected for the mathematics test represent a spread
across these situation types. In addition, items that can be regarded as
authentic are preferred.  That is, items should generally be mathematically
interesting and should reflect problems that could be encountered through a
person's day-to-day interactions with the world.
\end{enumerate}

In a typical Competency Class 1 problem, students were asked to read information
from a graph representing a physical relationship (speed and distance of a
car). Students needed to identify one specified feature of the graph (the
display of speed), to read directly from the graph a value that minimised the
feature, and then select the best match from given alternatives.

For a Competency Class 2 problem, students were given a mathematical model (in
the form of a diagram) and a written mathematical description of a real-world
object (a pyramid-shaped roof) and asked to calculate the area of the base. The
task required students to link a verbal description with an element of a
diagram, to recall the area formula for a square with given sides, and to
identify the required information in the diagram. Students then needed to carry
out a simple calculation to find the required area.

A Competency Class 3 task required students to identify an appropriate strategy
and method for estimating the area of an irregular and unfamiliar shape, and to
select and apply the appropriate mathematical tools in an unfamiliar context.
Students needed to choose a suitable shape or shapes with which to model the
irregular area, know and apply the appropriate formulae for the shapes they
used, work with scale, estimate length, and to carry out a computation involving
a few shapes.

PISA 2000 results:

To summarise data from responses to a collection of such items, a five-level
performance scale with an overall mean of 500 was created.  The scale was
created statistically using an Item Response Modelling approach to scaling
ordered outcome data.  Initially the overall scale was used to describe the
nature of performance by classifying the nations in terms of overall
performance, and thus to provide a frame of reference for international
comparisons.

For PISA 2000, the rank-order of countries showed that 15-year-olds in Japan
displayed the highest mean scores, but they could not be distinguished with
statistical significance from scores in Korea or in New Zealand.  Other
countries that also scored above the OECD average were Australia, Austria,
Belgium, Canada, Denmark, Finland, France, Iceland, Liechtenstein, Sweden,
Switzerland, and United Kingdom.  Overall, there was considerable within-country
variation.

Of more importance was the relationship of other variables such as student
motivation and engagement, gender, family background, and socioeconomic
background to performance in mathematics:
\begin{itemize}
\item In most countries, because most 15-year-olds considered mathematics
irrelevant to their future, only a small proportion considered mathematics worth
pursuing.

\item Lack of interest in mathematics was associated with poorer student performance.

\item Males on average performed better than females on mathematical literacy,
but the advantage disappeared when comparing low performers.

\item Higher parental education and more social and cultural communication among
parents and their children were associated with better student performance.

\item Living with only one parent was, on average, associated with lower student performance.

\item The socio-conomic composition of a school's student population was an even
stronger predictor of student performance than individual home background.
\end{itemize}

In summary, the PISA 2000 results provided an interesting initial look at how
15-year-olds responded to a set of items constructed to assess mathematical
literacy; the differences in mean performance across countries, and potentially
important correlates of such performance.

\section{The quantitative-qualitative debate and the case of a small-scale study}

\vskip-5mm \hspace{5mm}

Since FIMS in 1964, it seems that the same questions have repeatedly been asked
in large-scale studies, and that qualitative strategies are still not well
considered, although as a result of the criticism of FIMS, SIMS was
conceptualized as an in-depth-study of the curriculum. For the first time,
issues such as those related to student and teacher beliefs were
discussed. TIMSS added to SIMS studies, which aimed to explore the relationships
between the intended, implemented, and attained curriculum.

For example, in the study Survey of mathematics and science opportunities \cite{3},
based on observations in mathematics and science in several countries, it was
argued that there were typical patterns of instructional and learning activities
in each country, which seemed to stem from the interaction of curriculum and
pedagogy. It was assumed that students' learning experiences were moulded by
teachers who selected, prepared and taught the mathematical content in a variety
of instructional activities. In this respect, the researchers felt it necessary
to elicit information on teachers' background knowledge, their beliefs about
subject matter, and pedagogical beliefs. This view led the researchers to
explore teachers' instructional practices in detail. As the description of
teachers' practices through observations became the major focus, a major
reorientation in paradigm and methodology was inspired. The orientation of
conceptualization shifted from quantitative to qualitative differences between
countries. The Case study project \cite{13} was already designed at the beginning as
an ethnographic study, aiming to combine large-scale surveys and qualitative
methods. This in-depth studies of local situations intended to identify the
myriad of causal variables that were not recognised in large scale surveys and
to allow the development of hypotheses in order to interpret and explain many
data gathered in large scale studies.

Apart from these qualitatively oriented case studies accompanying TIMSS, there
exist many small scale international studies on mathematics education, and many
more are coming. As an example of such studies, we briefly describe one study
comparing English and German mathematics teaching \cite{5, 7}. This ethnographic
study was carried out at the beginning of the 1990s, using methods of
qualitative social sciences, mainly participating classroom observations. In
general, the study aimed at generating general knowledge, based on which
pedagogical phenomena might be interpreted and explained. Under a narrower
perspective, the study aimed to generate hypotheses on the differences between
teaching mathematics under the educational systems in England and in Germany.
For methodological reasons, however, the study could not make any ``lawlike"
statements; in contrast, the study referred to the approach of the ``ideal
typus" developed by Max Weber (Webersche Idealtypen) and described idealized
types of mathematics teaching reconstructed from the classroom observations in
England and Germany. That means that typical aspects of mathematics teaching
were reconstructed on the basis of the whole qualitative studies rather than on
an existing empirical case.

In brief, the study concluded that the following general approaches concerning
the understanding of mathematical theory were predominant in English and German
mathematics education.  In German mathematics teaching, a subject-oriented
understanding of mathematical theory prevailed, in contrast to the prevalence of
a pragmatic understanding of mathematical theory in England. Generally speaking,
mathematics teaching in Germany was characterized by its focus on the subject
structure of mathematics and on mathematical theory.  This meant that theory was
made explicit by means of rules and computations. In contrast, in England, the
understanding of theory could be called pragmatic---theory was applied
practically in an appropriate way.  These different basic approaches in England
and Germany to teaching mathematics were visible when looking at the differences
with regard to the following aspects.

The focus on theory when teaching mathematics in Germany implied a lesson
structure which went along with the subject structure of mathematics. Thus, in
the lessons, large units were complete in themselves. Mathematical theorems,
rules and formulae were therefore of high importance. That varied, though, with
the different kinds of schools of the three-track system.

The approach in England, the pragmatic understanding of theory, was apparent
from the curricular structure, which resembled a spiral. As a consequence,
smaller units were taught, and they were not necessarily connected with each
other.  Topics were quickly swapped, and at times different topics are worked on
at the same time.  Frequent repetitions of mathematical terms and methods that
had already been taught were a feature of this spiral-shaped approach.
Mathematical theorems, rules, and formulae (often called ``patterns") were of
low importance for the teaching of mathematics in England.

In Germany, proofs of mathematical statements were to a certain extent important
when teaching mathematics at the schools of the higher achievement level, but
they had only small or nearly no importance in schools of the intermediate or
lower achievement level. Proofs were considered important in order to visualise
the theoretical frame of mathematics, especially in the context of geometry.

In England, proofs were of low importance, both in selective as well as
in non\-selective schools.  Theorems, found by means of experiments,
were often only checked with examples, and proofs and checks with
examples were often not distinguished.

The status and role of proofs in Germany was studied more extensively in another
qualitative small-scale study by Knipping \cite{8}, comparing the proof of the
theorem of Pythagoras in French and German mathematics classes.

\section{Strengths and limitations of international studies in mathematics}

\vskip-5mm \hspace{5mm}

The strengths of comparative education can be seen in the
multidimensional aims described at the beginning such as describing
similarities and differences in educational phenomena in different
educational systems, estimating the effect of special variables on the
outcomes, based on input-output-models of education, identifying general
principles concerning educational effects. One major stream of
comparative work concerns itself with the interaction of educational and
political, social, or economic systems, for which the PISA study is an
example providing politicians with educational indicators, which aim to
steer educational systems. Another stream focuses on particular
pedagogical factors, for which the different case studies accompanying
TIMSS are an example. Comparisons of instructional methods, curricula,
teacher training, and their presumed outcomes (student behavior,
especially achievement) have long been at the heart of comparative work,
although the focus of attention has recently broadened. Alexander
\cite[p.~149]{1} describes this broadened view as follows:
\begin{quote}
I argue that educational activity which we call pedagogy---the purposive mix of
educational values and principles in action, of planning, content, strategy and
technique, of learning and assessment, and of relationships both instrumental
and affective---is a window on the culture of which it is a part, and on that
culture's underlying tensions and contradictions as well as its
publicly-declared educational policies and purposes.
\end{quote}

On the other hand, the limitations or even dangers of international comparative
studies are also now widely discussed among the researchers.

Kaiser \cite{6} described the following alternative approaches in comparative
education, which challenged established research traditions in comparative
education since the 1980s. The first challenge to the nation-state as the
exclusive research framework either looked at the world system---regional
variations, racial groups, classes, which are not bound to the nation---or did
microanalytic research focused on regional variation.  Proponents of the
analysis of regional variation argued that educational variance often was as
great, if not greater, across regions within a nation as it was across nations.

The second challenge to input-output models and total reliance on quantification
was based on assertions that education and school practices could not be reduced
solely to quantitative aspects---knowledge about these topics could only be
generated by qualitative research methods that focused on actual, lived
educational practices and processes.  A few approaches proposed
ethnomethodological techniques and related educational processes to broader
theories of school-society relations.

The third challenge questioned the dominance of structural functionalism in
comparative education---either how education functioned to maintain the social
fabric, or how it could be made to function (in the case of the Third World) to
develop a nation-state generally along Western models.  New approaches proposed
conflict theories, because most societies are plural societies characterised by
conflict, in which dominant groups seek to legitimise their control over the
state.

The fourth challenge, the emergence of new research concerns, involved new ways
of looking at educational institutions and their relation to society---studies
on the nature of knowledge transfer and its impact on the Third World, on school
knowledge, and on the internal workings of the school.

Concerning the international comparisons in mathematics education Kaiser \cite{6}
described four critical aspects, which should be considered in the coming
debate:
\begin{itemize}
\item On the methodological level, it is necessary to ask whether the main
concepts and ideas of the methodology used, such as the approach of the
probabilistic test theory, are adequate.  This model delivers highly general
results about ability levels.  It is, on the other hand, worth questioning
whether the necessary general conditions required by the model---existence of a
one-dimensional construct ``mathematical ability" or ``mathematical
literacy"---are fulfilled.

\item On a more subject-bound level, curricular issues have to be
discussed---the dependence on the results from the test items, the adequacy of
the test items concerning so many different curricula, the relation of
achievement results, and the opportunity to learn.

\item Under a general pedagogical perspective, there is a question of the
innovation potential of such studies.  What can we learn from descriptions of
mathematics teaching in totally different cultures, such as Germany and Japan,
with very different value systems, social conditions, and so on?

\item Under a political perspective, we need to ask whether the scientific
community is able to control how the results of the studies are used in
political debates. Consensus exists among researchers that the ranking is not
the main result of such studies---but the public debate concentrates mainly on
the ranking lists and bases proposals for consequences on the rank achieved.
How far are the results of such studies under the control of the researchers
involved, and how far can researchers influence the usage of the results?
\end{itemize}

Other branches of criticism emphasize the influence of multiple-choice
tests on mathematics education and on the promotion of mathematically
talented youngster: We have to distinguish between the mathematical
education for youngsters, who are highly interested and talented in
mathematics and intend to become a mathematician, engineer, etc. and
mathematics education for the others, who will not be necessarily
involved in vocations dealing with mathematics. The goals and the
methods of those two different branches of mathematical education are
quite different and should therefore use different methods of
comparison.

In many countries of the world, there is an elaborate system of
mathematical olympiads. They help to popularize mathematics, to
engage talented students into mathematics. The results of students
in major mathematical competitions (such as IMO) are sometimes
used to compare the levels of education in different countries.
The results of the team in IMO might be one indicator for the
level of mathematical education in a country. In addition, there
is another extraordinary international competition---the
Tournament of the Towns, ---which has a stronger focus on broader
groups of mathematically interested and talented students. It
offers an opportunity for students of many countries to compete
with other students in solving ambitious interesting mathematical
problems.

Coming back to the already mentioned large-scale international studies,
whose main part of the items are still multiple-choice items. Those
items are in many respects not adequate to compare the achievements of
mathematically talented students: One of the main ideas of classical
mathematical education looking for those students assumes that a student
is being taught not only to find solutions for different kinds of
problems, but to create different ways of solving problems, to create
arguments for the solution and communicate it orally and in written form
to an audience. Such a student is used to have a relatively large amount
of time per problem and he/she will have significant difficulties in
solving multiple-choice items under time pressure, where no
sophisticated argumentation is asked.

We have to consider that mathematics education in different countries is very
sensible to means of testing. Thus bringing multiple-choice tests into education
often leads to harmful changes in the whole education process. In Russia for
example, it has led to a degradation of an ambitious mathematics education in
many schools. It is an open question so far, whether mathematical olympiads or
other tournaments might be an appropriate indicator for the level of
mathematical education achieved by the students of a country, at least by the
mathematically talented, which could to a certain extent replace large-scale
international comparisons.

To summarize, comparative education is characterized by a wide diversity
of approaches, perspectives, and orientations, and this diversity of the
field seems to be one of its main strengths.

\label{lastpage}

\end{document}